\documentclass[12pt,pdftex]{amsart}
\usepackage[pdftex]{graphicx}

\usepackage{amssymb}

\newcounter{deficislo}[section]

\newcommand{\prop}[1]{\refstepcounter{deficislo}{\noindent \bf
Proposition \thesection.\thedeficislo.\  }{\it #1}}

\newcommand{\qz}{Q(z,\bar z,u)}
\newcommand{\g}{\gamma}

\newcommand{\ab}[1]{\vert z\vert^{#1}}
\newcommand{\cdva}{{\mathbb C^2}}

\newcommand{\te}{\theta}

\newcommand{\fz}{{F(z,\bar z,u)}}
\newcommand{\fzs}{{\textstyle{F^*(z^*,\bar z^*,u^*)}}}

\newcommand{ \al}{\alpha}
\newcommand{\be}{\beta}

\newcommand{\tf}{\tilde f}
\newcommand{\tg}{\tilde g}

\begin{document}
\title{
 Local symmetries of finite type hypersurfaces in $\mathbb C^2$}
\author {Martin Kol\'a\v r}
\address {Department of Mathematical Analysis, Masaryk University,
Janackovo nam. 2a, 662 95 Brno } \email {mkolar@math.muni.cz }
\maketitle

\begin{abstract}
The first part of this paper gives a complete description of local
automorphism groups for Levi degenerate hypersurfaces of finite
type in $\mathbb C^2$.
We also prove that, with the exception of  hypersurfaces of the form $v=\ab k$,
 local automorphisms are always determined by their 1-jets.
 Using this
result, in the second part we describe special normal forms
which by an additional normalization 
eliminate the nonlinear symmetries of the model and allow to decide effectively about
local equivalence of two hypersurfaces given in this normal
form.
\end{abstract}

http://front.math.ucdavis.edu/

\section{Introduction}
The main  aim of this paper is to give a complete description of
local automorphism groups for Levi degenerate hypersurfaces of
finite type in complex dimension two. The results rely mainly on
the construction of normal forms  given in \cite{Kol}.

The problem of describing local  symmetries for real hypersurfaces
in two dimensional complex space is closely related to the local
equivalence problem. This connection appears already in the
foundational work  of H. Poincar\'e \cite{Po}.
 A natural approach
to the problem, whose germs can be also found in \cite{Po},  is to
analyze directly the action of the group  of local biholomophic
transformations at the given point. For Levi nondegenerate
hypersurfaces this analysis was completed in the beautiful
construction of S. S. Chern and J. Moser \cite{CM}.

In recent years the same approach was applied to various classes
of Levi degenerate hypersurfaces  (see \cite {BB, BE, E, ELZ2, J, S, W}). In
\cite{Kol}, normal forms are constructed for general finite type
hypersurfaces in dimension two. Since the normal forms are on
the level of formal power series (neither convergence nor
divergence  has been proved),
  in order to solve Poincar\'e's local equivalence problem
one has to combine the construction with the result of M. S.
Baouendi, P. Ebenfelt and L. P. Rothschild \cite{BER1} on
convergence of formal equivalences between finite type
hypersurfaces. In fact, there are three different normal forms
defined in \cite{Kol}, depending on the form of the model (see
Section 2 for more details).

 For Levi nondegenerate
hypersurfaces, Chern-Moser's construction of normal forms,  in
general dimension, gives already substantial  information about
local symmetries, but not complete information. It shows that the
local automorphism  group of any hypersurface is a subgroup of the
group of local symmetries of the model hyperquadric. It also
implies that local automorphisms are detemined by their 2-jets.

 These results are strengthened and
completed in an important way by the theorems of V. K. Beloshapka,
N. G. Kruzhlin and A. V. Loboda (\cite{B, KL, V}).  Local
automorphisms are in fact determined by their 1-jets, whenever the
hypersurface is different from the model hyperquadric.  Moreover,
in the strongly pseudoconvex  case there exist local holomorphic
coordinates in which all automorphisms are linear. Results of V. V. Ezhov
(\cite{Ez}) show that the last property no longer holds in general for hypersurfaces 
with mixed signature. 

The group of local automorphisms of a hypersurface (i.e.
automorphisms which fix the given point) is usually called the
stability group. The problems of finite jet determination and
estimation of the dimension of the stability group on Levi
degenerate hypersurfaces have been intensively studied in the last
decade (see \cite {EHZ, ELZ, ELZ2, EL, Z},  and the  survey article
\cite{BER2} for further references). In dimension two, one of the
most important results states that uniform finite determination,
which holds on finite type hypersurfaces, actually fails for
points of infinite type (see \cite {Kow,Z}). More precisely, for
any integer $k$ there is an infinite type, non Levi flat
hypersurface whose local automorphisms are not determined by their
$k$-jets. On the other hand, D. Zaitsev formulated recently a conjecture
 that in dimension two,
jets of order higher than one are needed only for hypersurfaces which are biholomorphic 
to the ball at generic points.
Proposition 3.1 below confirms this conjecture in the finite type case.

In Section 2 we introduce notation and review the needed
ingredients of the normal form construction from \cite {Kol}. In
Section 3 we 
consider hypersurfaces whose models are 
 the higher type analogs of spheres, given by
$v = \ab k$. We prove that, except for the model hypersurfaces themselves,
 local automorphisms are always determined by 1-jets,
and are linear in normal coordinates (it should be stressed that
since the convergence of the normal forms  has not been proved,
these normal coordinates are a priori only formal). In Section 4
we use this result to give a complete classification of local
automorphism  groups for Levi degenerate hypersurfaces of finite type.
The remaining analysis for non-spherical models is
straightforward, as the symmetries  of such models are themselves
linear. As one consequence, the result gives a complete description
of hypersurfaces with finite stability group.  In Section 5  we define special normal forms for
hypersurfaces whose  models are given by $v = \ab k$, which allow
to decide effectively about local equivalence of two hypersurfaces
 put into this  normal form.

I would like to thank the  organisers of the 
2005 International Conference on Several Complex Variables in Hefei 
 for their invitation
and hospitality. I am also grateful to Peter Ebenfelt, Bernhard Lamel 
and Dmitri Zaitsev for helpful discussions. 

During the work on this paper I learned from M. Eastwood about the
work of Valerii  Beloshapka and Vladimir Ezhov \cite{BE}. In this paper
normal forms similar to those of \cite{Kol} are described (six
 different cases are considered).

\section{Notation and preliminaries  }

We will consider a real analytic hypersurface $M\subseteq \mathbb
C^2$ in a neighborhood of a point $p \in M$ in which the Levi form
degenerates. The point will be  assumed to be of finite type $k$
in the sense of Kohn (\cite{K}).

 For local description of $M$ near $p$ we will use local
holomorphic coordinates $(z,w)$, where  $z = x + iy,\  w = u +
iv$, such that   $p = 0$ and the hyperplane $\{ v=0 \}$ is tangent
to $M$ at $p$. $M$ is then described near $p$ as the graph of a
uniquely determined real valued function
\begin{equation}
v= F(x,y,u).\label{2.1}
\end{equation}
 We will consider the Taylor expansion of
$F$ expressed in terms of $(z, \bar z, u)$:
\begin{equation}\fz = \sum_{i+j+m\ge 2} a_{ijm}
 z^i \bar z^j u^m,\label{2.3}\end{equation}
where $a_{ijm}=\overline{a_{jim}}$.

Further we denote
\begin{equation}Z_{ij}(u) = \sum_{m} a_{ijm} u^m,\end{equation}
hence
\begin{equation}F(z, \bar z, u) = \sum_{i,j}Z_{ij}(u)z^i\bar z^j.\end{equation}
 We will analyze the effect of a
  holomorphic transformation
\begin{equation}  z^*=f(z,w), \ \ \ \ \   w^*=g(z,w),  \label{2.4} \end{equation}
on the defining equation of $M$. Here $f$ and $g$ are represented
by power series
\begin{equation} f(z,w)=\sum_{i,j =0}^{\infty} f_{ij} z^i w^j, \ \ \ \ \   g(z,w)= \sum_{i,j=0}^{\infty} g_{ij}
 z^i w^j.
 \label{2.5} \end{equation}

We require that such a transformation  preserves the form  given
by (\ref{2.1}), which means that    the origin has to be mapped to
itself and the hyperplane  $v^* = 0$ has to be tangent to $M$ at
$p$ in the coordinates $(z^*, w^*)$.
 This will be satisfied if and only if
\begin{equation}f=0,\ \ g=0, \ \ g_z=0, \ \ Im\ g_w=0\ \ \  \text{at}\ \ \
z=w=0.\label{2.6}\end{equation}  Only such transformations will be
considered in the sequel.

 We denote by $F^*$ the function describing $M$ in coordinates
$(z^*, w^*)$, and write
\begin{equation}\fzs = \sum_{i+j+m\ge 2} a_{ijm}^* (z^*)^i (\bar z^*)^j (u^*)^m,\label{2.3}\end{equation}
where again $a^*_{ijm}=\overline{a^*_{jim}}$.
$F^*$ is related to $F$ and $(f,g)$  by the following
transformation formula:
\begin{equation}  F^*(f(z, u+iF),\; \overline{ f(z,u+ i F)},\;
Re\; g(u + iF))
  = Im\ g(z,u+iF),
\label{2.7} \end{equation} where the argument of $F$ is $(z, \bar
z , u).$
This formula gives an equality of two power series in $z, \bar z, u$,
and by comparing coefficients of various monomials we may obtain 
 explicit relations between $F^*$ and $F, f, g$.

A natural tool which simplifies the use of this formula is provided by weighted
coordinates.
We give weight one to $z$ and $\bar z$ and weight $k$ to $u$ and $w$.

Recall that $p \in M$ is a point of finite type $k$ in the sense of
Kohn if and only if there exist local holomorphic coordinates such
that $M$ is described by
\begin{equation}v = \sum_{j=1}^{k-1}
a_j z^j\bar z^{k-j} + o(\ab k, u ),\ \label{2.9}\end{equation}
where the leading term  is a nonzero real valued homogeneous
polynomial of degree $k$, with $a_j \in \mathbb C$ and $a_j =
\overline{a_{k-j}}.$

The model hypersurface $M_D$ to $M$ at $p$ is defined using the leading
homogeneous term,
\begin{equation}M_D = \{(z,w) \in \cdva\ | \  v  = \sum_{j=1}^{k-1}
a_j z^j\bar z^{k-j} \}. \label{2.10}\end{equation}
 In particular, when the leading term is equal to  $ \ab k$, we will write
\begin{equation}O_k = \{(z,w) \in \Bbb C^2 \ \vert \   v = \ab k
\}. \end{equation}

Now we recall two basic integer invariants used in the normal form
construction in \cite{Kol}. The first one, denoted by $l$, is  the
essential type of the model hypersurface to $M$ at $p$. It can be
described as the lowest index in (11) for which $a_l \neq 0$.
It satisfies $1 \leq l \leq \frac k2.$ 

The second invariant is defined when $l < \frac k2$ as follows. Let $l=m_0 <m_1 <\dots
<m_s <\frac k2$ be
 the indices in (11) for which $a_{m_i}\neq 0$. The invariant,
 which we will denote
 by $\kappa$,
 is the greatest common divisor of the numbers
$\ k-2m_0, k-2m_1,  \dots, k-2m_s$.

While for $l < \frac k2$ the stability group of $M_D$ is one dimensional, 
the stability  group of $O_k$ has dimension three. Its
elements are of the form $(\tf, \tg)$, where
\begin{equation}\tf(\delta, \mu, \theta; z,w)  =  \frac{ \delta e^{i\theta} z}{(1 + \mu w)^{\frac1l}},
\ \ \ \ \tg(\delta, \mu, \theta; z,w)=
  \frac{ \delta^k w}{1 + \mu w},\label{4.1}\end{equation}
                           with $\delta > 0,$ and $\te, \mu  \in \Bbb R$.
 We will use their Taylor expansion
\begin{equation}\tf(\delta, \mu, \theta; z,w) =
\delta e^{i\theta}( z - \frac{\mu}l zw + \dots \ )\end{equation}
and
\begin{equation}\tg (\delta, \mu, \theta; z,w) = \delta^k
( w - \mu w^2 + \dots \ ).\end{equation}
 In \cite{Kol}, Proposition 4.2,  we proved that if the model hypersurface to $M$ at $p$ is $O_k$,
 there exists a unique formal
 transformation satisfying normalization conditions (7) and
\begin{equation} f_z= Re\; g_w = 1, \ \ \ Re \; g_{ww} = 0  \ \ \ at\  z=w=0,\label{3.8}\end{equation}
which takes the defining equation for $M$ into normal form, where
the normal form conditions are

$$
\begin{aligned}
Z_{j0}&= 0, \ \ \ \ \ j=0,1,\dots,  \\
Z_{l,l+j}&= 0, \ \ \ \ \ j=0,1,\dots, \\
Z_{2l, 2l}& = 0,\\
Z_{3l, 3l}& = 0,   \\
Z_{2l, 2l-1}& = 0. \ \end{aligned}
$$

\section{Linearity of local automorphisms}

In this section we prove a result analogous to the result of
 \cite{B, KL} for nondegenerate hypersurfaces.
The stability group of $M$ at $p$ will be denoted by $H_M$.
In \cite{Kol} we proved that if $l < \frac k2$, then $dim\;  H_M \leq
1$, and all local automorphisms are determined by their 1-jets.
Here we consider the case when $l = \frac k2$.

\prop {\it If $M$ is not equivalent to $O_k$, then $dim\;  H_M \leq
1$. Moreover, all local automorphisms expressed in normal
coordinates are linear.}

 {\it proof:}
We assume
 that the model is $O_k$, but $M$ is not equivalent to $O_k$.
  Let us consider normal coordinates for $M$ and
separate the first two  leading terms in the Taylor expansion of
$F$,
\begin{equation} \fz = \ab k + \qz + o_{wt}(p), \end{equation}
where $Q$ is a nonzero weighted homogeneous real valued polynomial
of weight $ p > k, $
\begin{equation} \qz = \sum_{\al + \be  + k \g = p} a_{\al \be \g}z^{\al} \bar
z^{\be} u^{\g}, \end{equation} and $o_{wt}(p)$ denotes terms which
are of weight greater then $p$. We define the index $(\al_0,
\be_0, \g_0)$ to be the smallest one in inverse lexicographic
ordering (the last components are compared first, then the second
ones) for which $a_{\al_0, \be_0, \g_0} \neq 0$.

 Let $(f,g)$ be a local automorphism of $M$,
i.e. a transformation which preserves $F$. Its general form 
is
\begin{equation}
\begin{aligned}
f(z,w) &=\delta e^{i\theta} z + \text{terms of weight} \ge 2,
\\
 g(z,w) &= \delta^k w + \text{terms of weight}
 \ge k+1.
\end{aligned}
\end{equation}
 We will call the
numbers $\delta, \; \theta \; $ and $ \mu = Re\; g_{ww}$
 the initial data of the automorphism, and  consider
simultaneously  $M$ with the automorphism $(f,g)$ and the model $O_k$ with the
automorphism $(\tf, \tg)$ having  the same initial data as
$(f,g)$. We will use (9) to compare the coefficients of $(f,g)$ and 
$(\tf, \tg)$ (see [18] for a detailed description of the use of the transformation 
rule (9)).  

 First we will show that $f$ and $g$ may
be replaced by $\tf$ and $\tg$ when considering terms of weight
less or equal to $p+k$ in (\ref{2.7}). More precisely,
\begin{equation}f(z,w) =  \tf(z,w) + o_{wt}(p+1),\ \ \ \ \
 g(z,w) = \tg(z,w) + o_{wt}(p+k).\end{equation}
This is done in two steps. First, since $Q$ has weight p, 
all
equations obtained from (\ref{2.7}) for coefficients of monomials  up to weight $p-1$ are the
same as those for $O_k$ and $(\tf$, $\tg)$. Hence $f$ is equal to
$\tf$ modulo $ o_{wt}(p-k)$ and $\tg$ equal to $g$ modulo $
o_{wt}(p-1)$. For terms of weight $p$, $Q$ enters (\ref{2.7}) only  via
the linear part of $(f,g)$, as
$Q(\delta e^{i\theta} z, \delta e^{-i\theta}  \bar z, \delta^k u)$.
 Since $Q$ (and in particular  $a_{\al_0,
\be_0, \g_0}$)  has to be preserved, we
obtain immediately that 
 $$\delta = 1, \ \ \ \ \ e^{i(\al_0 - \be_0 )
\theta} = 1.$$
 For terms of weight $p+1, p+2, \dots, p+k$, $Q$
enters (\ref{2.7}) only through  the initial data
$Re g_{ww}$, and   the coefficients   $f_{20}, \dots f_{k0}$
in $f$ and $g_{11}, \dots  g_{k1}$ in $g$.
 But we
already know these coefficients to be the same as in $(\tf, \tg)$, namely 
 zero (if $ k > p-k$ 
we use  an obvious step by step argument).
Since by the result of \cite{Kol}  a local automorphism is uniquely determined 
by its initial data, it follows that $\tf$ has to agree with $f$
modulo terms of weight greater than $p+1$ and $\tg$ has to agree with $g$ modulo
terms of weight greater than $p+k$. 
 This proves the
claim.

Now we consider all terms of weight $k+1, \dots, k+p$ in the
transformation formula (\ref{2.7}) . On the right hand side,
using $ g(z,w) = w - \mu w^2 + \dots$
 we
have

\begin{equation} \begin{aligned} &Im \; g(z, u + iF) = 
 F - Im\;  \mu (u + i(\ab k + Q +
o_{wt}(p)))^2 + J_1 + \\ +& \; o_{wt}(k+p)  = F + 2\mu u \ab k - 2 \mu u Q + J_1 +
o_{wt}(k+p),
\end{aligned}\end{equation}
where  $J_1$  denotes  terms of weight $\leq k+p$ which come
only from $\ab k$, in other words terms which appear in the corresponding 
expansion for $O_k$ and $(\tf, \tg)$
  (which we will not need to write down explicitely).
 On the left,
 we get from the leading term
\begin{equation}
\vert f(z, u + i(\ab k + Q + o_{wt}(p)))\vert^k = 
\vert z - \frac{\mu}l z( u + i(\ab k + Q + o_{wt}(p)))+ o_{wt}(2k)\vert^k  
\end{equation}
which gives
\begin{equation}
 \ab k
- 2\mu Im \; \ab k Q + J_2 + o_{wt}(p+k)
=
 \ab k
 + J_2 + o_{wt}(p+k)
,\end{equation}
where again $J_2$  denotes all terms of weight $\leq k+p$ which
come only from $\ab k$.
 From the second term in $F=F^*$ we get 
\begin{equation}
\begin{aligned}
 &Q(f,\bar f, Re\ g) = Q(
e^{i\theta} z-e^{i\theta}\frac{\mu}l z(u+ i \ab k)  + o_{wt}(k+1)), \bar{
\ \ }, 
\\& u - Re\; \mu (u +
i(\ab k + o_{wt}(k)))^2+ o_{wt}(2k)). \end{aligned}
\end{equation}

By the same argument as we used before for $Q$,
since   $f_{20}, \dots f_{k0}$
 and $g_{11}, \dots  g_{k1}$ vanish,
 terms of weight greater than 
$p$ and less or equal to $p+k$  in $F^*$ influence  (9) only via the
 linear part of $(f,g)$.
Multiplying out and
taking into account that terms coming only from 
$\ab k$ have to  eliminate each other,
we calculate the coefficients of $z^{\al_0} \bar z^{\be_0}
u^{\g_0+1}$ in (9). We get
\begin{equation}
 a_{\al_0,
\be_0, \g_0 +1} - a_{\al_0, \be_0, \g_0}( \frac1l \mu (\al_0 + \be_0 +l  \g_0))
 = - 2\mu a_{\al_0, \be_0, \g_0 } +
a_{\al_0, \be_0, \g_0 +1}.
\end{equation}
It will hold if and only if
\begin{equation}
a_{\al_0,  \be_0,  \g_0 }
(2\mu - \frac1l \mu ( \al_0 + \be_0 + l \g_0 )) = 0,\end{equation}
 hence
\begin{equation}\al_0 + \be_0 + l \g_0 = k\end{equation}
(recall that $ l = \frac k2$).
It follows that either $\mu = 0$, or  $\g_0 = 1$ and $\al_0 + \be_0 = l$.
If $\g_0 = 1$ we
consider the coefficients of $z^{\al_0 + k}\bar z^{\be_0 + k}$. From 
the formulas above  
we get
\begin{equation}a_{\al_0 + k,\be_0 + k,0} + \mu a_{\al_0, \be_0, 1} = a_{\al_0 + k,\be_0 +
k,0},\end{equation}  and so $\mu = 0$. Hence there is no
automorphism of $M$ with $\mu \neq 0$, and we proved that every
local automorphism in normal coordinates is linear. To prove that $dim \
H_M \leq 1$, it is enough to realize that linear automorphisms act
on each term in $F$ individually, and that dilations preserve only
the homogeneous model.

\section{Classification of local symmetries}

The result of Proposition 3.1 can be used to obtain a complete
classification of local automorphism groups.

\prop  {\it For a given hypersurface $M$ exactly one of the following
possibilities occurs:} \begin{enumerate} {\it \item $H_M$ has real
dimension three. This happens  if and only if $M$ is equivalent to
$O_k$.

 \item $H_M$ is isomorphic to
$\Bbb R^+ \oplus {\Bbb Z}_m$  This happens if and only if $M$ is a
model hypersurface with $l<\frac k2$, and $m = \kappa $ when $k$ is even 
or  $m = 2\kappa $ when $k$ is odd.

\item $H_M$ is isomorphic to $S^1$.  This happens if and only if
$M$ is weakly spherical, i.e. the defining equation in normal
coordinates has form
$$ v = G(\ab 2, u).$$
\item $H_M$ is finite, isomorphic to $ {\Bbb Z}_n$ for some $n \in \Bbb N$.
}
\end{enumerate}
 Note that the last
case includes the trivial symmetry group.

{\it proof:}   By Proposition 3.1, if $l = \frac k2$ and $M$ is
 not equivalent to $O_k$,
then the only transformations which may preserve it in normal
coordinates  are the decoupled linear transformations
\begin{equation} w^* = \delta^k w , \ \ \ \ \ \ \ z^* = \delta e^{i\theta} z, \ \label{3.3}\end{equation}
which act on each term in the expansion of $F$ individually in an
obvious way. Each such transformation can be uniquely factored
into the composition of a rotation in the $z$ variable and a
weighted dilation.  If the rotation $z^* = e^{i\theta} z$
preserves all terms in $F$ for every $\theta$, then each must have
form $a_{m_1, m_1, m_2}\ab {2m_1} u^{m_2}$, corresponding to the third case.
Further,  if
for one particular $\theta$ the  rotation preserves a term $a_{\al \be \g} z^{\al} \bar z^{\be} u^{\g}
 =  \ab
{2\al}Re \; a_{\al \be \g} z^{\be - \al} u^{\gamma}$,
 where $\al \leq \be$, then
$e^{i\theta} $
is a $(\be - \al)-th$  root of unity, and we are in cases (2) and (4).
On the other hand, in all cases, a weighted dilation can preserve only terms
which are weighted homogeneous of weight $k$. In this case $M$
has to be a model. Further, for $l < \frac k2$ 
it follows from \cite{Kol} that $H_M$ is a subgroup of 
$\Bbb R^+ \oplus {\Bbb Z}_m$, 
with the claimed relation between $m$ and $\kappa$.
Hence, if $M$ is not a model, $H_M$ has to be a subgroup of ${\Bbb Z}_m$,
i.e. it is isomorphic to ${\Bbb Z}_n$ for some $n \in \Bbb N$. 

\section{Special normal forms}

In this part we use the calculations from Section 2 to obtain a
normal form which can be used effectively to decide about local
equivalence of two hypersurfaces given in normal form. The main
difficulty in applying Chern-Moser's normal form for that purpose
is not the dimension of the symmetry group, but rather the fact
that the group does not act on the defining equation directly.
Application of an element of the group can lead to an equation not
in normal form. To obtain the group action on normal forms one has
to perform the transformation back into normal form.
~The same situation occurs for the normal forms in \cite{Kol}, in
the case when $l = \frac{k}2$ and the model is $O_k$.

We will speak about a special normal form if it practically allows
to decide about equivalence or non-equivalence of two
hypersurfaces which are put into this normal form. More precisely,
in the two dimensional case this means that only the explicit
action of a (decoupled) linear transformation is to be considered.
In the nondegenerate case such a
normal form is described in \cite{CM} for non-umbilical points in
$\mathbb C^2$ and in \cite{W} for  higher dimensions.

 We rewrite $Q$ in the
form

$$ \qz = \sum_{\al + \be  + k \g = p} \ab {2\al} Re \;
a_{\al \be \g}  z^{\be - \al}   u^{\g}. $$
where the sum is taken
over multiindices with $\alpha \leq \be$.
 Recall that  $({\al_0, \be_0, \g_0})$ is the
smallest index in (11) in inverse lexicographic ordering for which $a_{\al_0, \be_0, \g_0}$
 is different from
zero.

First we normalize the linear part of a transformation into normal
form by requiring that
$$a_{\al_0, \be_0, \g_0} = 1.$$
This condition provides a partial normalization or a complete one,
depending on the form of $Q$. In all cases it normalizes fully the
dilation part, while $\theta $ is left free if $Q$ consists of a
single term of the form $\ab m u^s$. In the second step we
normalize the nonlinear part. As in the proof of Proposition 3.1
we have to consider two cases. If 
\begin{equation}\al_0 + \be_0 + l \g_0 \neq k,\end{equation}
we normalize  by asking that

\begin{equation} Re \; a_{\al_0, \be_0, \g_0 +1}\ = 0.\end{equation}
In the second case, when $\g_0 = 1$ and $\al + \be = l$ we normalize
 by requiring  that 
\begin{equation} Re \;a_{\al_0 + k,\be_0 + k,0} = 0.\end{equation} 
By the calculation in the proof of Proposition 3.1, 
these conditions determines uniquely the parameter $\mu$.

  Thus verifying local equivalence of two such hypersurfaces is reduced 
 to the straightforward action of linear
transformations.

\vspace{5mm}
\begin{center}
{\textsc{Acknowledgement:}} \end{center}
 {Supported by
 a grant of the GA \v CR no. 201/05/2117}
\end{document}